\documentclass[12pt]{article}

\usepackage{amsmath}
\usepackage{theorem}
\usepackage{amssymb}
\usepackage{latexsym}
\usepackage{a4wide}
\usepackage{longtable}
\usepackage{epic}
\usepackage[all]{xy}

\newtheorem{thm}{Theorem}
\newtheorem{lem}[thm]{Lemma}
\newtheorem{prop}[thm]{Proposition}
\newtheorem{cor}[thm]{Corollary}

{\theorembodyfont{\rmfamily} }
{\theorembodyfont{\rmfamily} }
\newenvironment{rem}{\noindent{\bf Remark.}}{\newline}
\newenvironment{pf}{\noindent{\bf Proof.}}{\hbox{}\hfill $\Box$}
{\theorembodyfont{\rmfamily}}

\newcommand{\C}{\mathbb{C}}
\newcommand{\R}{\mathbb{R}}
\newcommand{\Q}{\mathbb{Q}}
\newcommand{\Z}{\mathbb{Z}}

\newcommand{\GL}{\mathrm{\mathop{GL}}}

\newcommand{\Tr}{\mathrm{\mathop{Tr}}}

\newcommand{\ad}{\mathrm{\mathop{ad}}}
\newcommand{\Stab}{\mathrm{\mathop{Stab}}}

\newcommand{\gl}{\mathfrak{\mathop{gl}}}
\newcommand{\mf}[1]{\mathfrak{#1}}
\newcommand{\ssl}{\mathfrak{\mathop{sl}}}
\newcommand{\g}{\mathfrak{g}}

\newcommand{\hh}{\mathfrak{h}}

\begin{document}

\title{Constructing semisimple subalgebras of semisimple Lie algebras}
\author{Willem A. de Graaf\\
Dipartimento di Matematica\\
Universit\`{a} di Trento\\
Italy}
\date{}
\maketitle

\begin{abstract}
Algorithms are described that help with obtaining a classification 
of the semisimple subalgebras of a given semisimple Lie algebra, up to
linear equivalence. The algorithms have been used to obtain classifications
of the semisimple subalgebras of the simple Lie algebras of ranks $\leq 8$.
These have been made available as a database inside the {\sf SLA} package
of {\sf GAP}4. The subalgebras in this database are explicitly given, as
well as the inclusion relations among them.
\end{abstract}

\section{Introduction}

It is an extensively studied problem to classify the semisimple
subalgebras of a complex semisimple Lie algebra $\g$, up to an equivalence
relation. The most natural equivalence relation for this is the one of
conjugacy by the inner automorphism group $G$. 
Two subalgebras of $\g$ are simply called 
{\em equivalent} if they are conjugate under
$G$. In \cite{dyn}, Dynkin also considered the notion of linear equivalence:
two subalgebras $\g_1,\g_2\subset\g$ are said to be {\em linearly equivalent}
if for every representation $\rho : \g \to \gl(V)$ the subalgebras 
$\rho(\g_1)$, $\rho(\g_2)$ of $\gl(V)$ are conjugate under $\GL(V)$. 

A subalgebra of $\g$ is called {\em regular} if it is normalised by a 
Cartan subalgebra of $\g$. Semisimple subalgebras of this kind correspond to 
root subsystems of the root system of $\g$. An $S$-{\em subalgebra} is a
subalgebra which is not contained in a regular subalgebra.   

In \cite{dyn0} Dynkin classified the maximal semisimple $S$-subalgebras 
of the Lie algebras of classical type, upto equivalence. 
More precisely, \cite{dyn0} contains a description, or a procedure, by which 
it is possible for a given Lie algebra of classical type to find its maximal 
semisimple $S$-subalgebras.

Dynkin treated the Lie algebras of exceptional type in \cite{dyn}. The
main results of this paper are
\begin{itemize}
\item an algorithm to classify the regular subalgebras of a semisimple
Lie algebra, up to equivalence,
\item a classification of the semisimple $S$-subalgebras, up to equivalence, 
of the Lie algebras of exceptional type,
\item a classification of the simple subalgebras, up to linear equivalence, 
of the Lie algebras of exceptional type.
\end{itemize}

Lorente and Gruber (\cite{logru}) applied Dynkin's methods to obtain explicit 
lists of semisimple subalgebras of the simple Lie algebras of classical type.
More in particular, they obtained lists of the regular subalgebras and of
the $S$-subalgebras of the Lie algebras of classical type of ranks $\leq 6$.

Recently Minchenko (\cite{minchenko}) has revisited Dynkin's classification
of the simple subalgebras of the simple Lie algebras of exceptional type.
He corrected several small mistakes (most notably he found two extra simple
subalgebras in the Lie algebra of type $E_8$). Secondly he found the 
classification of the simple subalgebras up to equivalence. Thirdly, he has
computed a lot of additional data (such as the normalisers of the subalgebras
in $G$). 

One motivation for studying semisimple subalgebras of semisimple
Lie algebras comes from theoretical physics. In models like the
vibron model and the interacting boson model (cf. \cite{iaar}) 
chains of subalgebras are
used. For applications of this kind the subalgebras need to be explicitly
given, i.e., for each equivalence class a representative needs to 
be given by a basis. Furthermore, methods are needed to obtain the 
inclusion relations between the subalgebras (more precisely: to decide
whether two given classes have representatives such that one is contained
in the other).
The classifications present in the literature do not appear to immediately
give this. For example, in \cite{dyn}, only the $S$-subalgebras
are explicitly constructed. And only the simple subalgebras are
listed, and not the semisimple ones (with the exception of the 
$S$-subalgebras). Finally no inclusion relations are given (again with the
exception of the $S$-subalgebras). 

The aim of this paper is to describe algorithms, and report on the results
obtained with their implementation, that help with obtaining a classification 
of the semisimple subalgebras of a given semisimple Lie algebra, up to
linear equivalence. Furthermore,
the subalgebras are explicitly constructed, as well as the inclusion 
relations among them. Here we say that the algorithms ``help'' to obtain
a classification as one step in the algorithms (the construction of the
subalgebras) is not entirely algorithmic - occasionally some human intervention
is needed for that.

Equivalence implies linear equivalence, but the converse is not always
true. However, if $\g$ is of type $A_n$, $B_n$, $C_n$, $F_4$, $G_2$ then
the two concepts coincide (cf. \cite{minchenko}, Theorem 3). 
In the remaining types there are some exceptions
and they are explicitly described (\cite{dyn}, \cite{minchenko}).
Hence it is straightforward to obtain the classification of the semisimple
subalgebras up to equivalence from the list of semisimple
subalgebras up to linear equivalence.
One of the main advantages of linear equivalence as opposed to equivalence
is that we have a method for deciding it (see Section \ref{p:6}). For these
reasons in this paper we focus exclusively on linear equivalence. 

By considering embeddings of Lie algebras in $\g$, rather than subalgebras
of $\g$ we get a slightly different perspective on the problem.
Also for embeddings we have the notions of equivalence and linear equivalence.
 Let $\tilde\g$ be a semisimple 
Lie algebra, and $\varphi_1,\varphi_2 : \tilde\g\hookrightarrow \g$ injective
homomorphisms. They are said to be equivalent if there is a $\sigma\in G$
with $\varphi_1 = \sigma\varphi_2$. They are said to be linearly equivalent if
for each representation $\rho : \g \to \gl(V)$ the induced representations
$\rho\varphi_1$, $\rho\varphi_2$ of $\tilde\g$ are equivalent.
Let $\tilde \g$ be a semisimple Lie algebra, and $\g' \subset \g$ a
subalgebra isomorphic to $\tilde\g$. There can be several
non-equivalent embeddings $\tilde \g \hookrightarrow \g'$. This is only
possible if $\tilde \g$ has outer automorphisms. From a classification
of subalgebras up to linear equivalence it is straightforward to get
all embeddings up to linear equivalence. For this reason
we concentrate on constructing subalgebras, rather than embeddings.

One approach to the problem is to start from the existing classifications in 
the literature. One could
take the maximal $S$-subalgebras constructed by Dynkin, along with the
regular subalgebras, and by successively
constructing their subalgebras get the entire list of subalgebras. This,
however, would not confirm, or correct, the existing classifcation. Moreover,
if the list of maximal subalgebras has an error, then this will lead to
many errors in the resulting classification. (And it appears that this can,
for example, easily happen in type $D_{2n}$, see below.) For 
these reasons the approach taken here aims at obtaining the classification
from scratch. This has the added advantage that the classifications in the 
literature and the new ones can validate each other. In particular, if
both are the same then this constitutes a good argument for their correctness.

The main idea used here to classify subalgebras is to start with the ones
of smallest rank. The subalgebras of rank 1 are well-known from the 
classification of the nilpotent $G$-orbits in $\g$. Secondly we construct
the subalgebras of higher rank as a kind of extension of algebras of lower
rank. This way we ``climb our way up''. So, in a sense, it is the reverse
approach to starting with the maximal subalgebras.

The algorithms described in this paper have been implemented in the
language of the computer algebra system {\sf GAP}4 (\cite{gap4}), using the
package {\sf SLA} (\cite{sla}). The main result that 
has been obtained using this implementation is a database of all
semisimple subalgebras of the simple Lie algebras of ranks $\leq 8$.
This database is also contained in the package {\sf SLA}. It
also contains all inclusion relations between the linear equivalence classes.
It is complemented by a function for computing the semisimple subalgebras
of a semisimple, non-simple,
Lie algebra. In Table \ref{tab:1} we show some statistics
relative to the simple Lie algebras of ranks $7$, $8$. The table contains
the number of (linear equivalence classes of) subalgebras, and the number of 
their isomorphism types. 

There is also the question of the field of definition. The simple Lie
algebras are given by a multiplication table relative to a Chevalley
basis. The subalgebras are given by a basis. However, not all linear
equivalence classes of subalgebras have a representative with a basis
with coefficients in $\Q$ (with respect to the given Chevalley
basis). Our results show that for all semisimple subalgebras of the
simple Lie algebras of ranks $\leq 8$ there exists a field extension $F$
of degree $\leq 2$ of $\Q$, such that the subalgebra can be given by a basis
with coefficients in $F$. The last column of Table \ref{tab:1} gives a field
extension $F$ of $\Q$ such that all semisimple subalgebras can be given by 
a basis with coefficients in $F$. Here we remark that it is by no means clear
that these are the smallest possible fields (except, of course, when the field
is $\Q$). We have made an effort to keep the fields small; but the problem of
finding the absolute smallest field is a difficult one which we do not solve
here.

\begin{table}[htb]\label{tab:1}
\begin{center}
\begin{tabular}{|l|r|r|l|}
\hline
& \# subalgebras & \# types & field of definition\\
\hline
$A_7$ & 131 & 32 & $\Q$\\
$B_7$ & 849 & 95 & $\Q(\sqrt{-1},\sqrt{-2},\sqrt{-3})$ \\
$C_7$ & 822 & 76 & $\Q$\\
$D_7$ & 511 & 72 & $\Q(\sqrt{-1},\sqrt{-3},\sqrt{-5})$ \\
$E_7$ & 501 & 76 & $\Q(\sqrt{-3})$ \\
$A_8$ & 232 & 46 & $\Q$ \\
$B_8$ & 2186 & 165 & $\Q(\sqrt{-1},\sqrt{-3})$ \\
$C_8$ & 2127 & 126 & $\Q(\alpha)$, $\alpha^2-\alpha-1=0$\\
$D_8$ & 1664 & 127 & $\Q(\sqrt{-1},\sqrt{-3})$\\
$E_8$ & 1183 & 155 & $\Q(\sqrt{-1},\sqrt{-3},\beta)$, $\beta^2+\beta+2=0$\\
\hline
\end{tabular}\caption{Semisimple subalgebras of the simple Lie algebras of rank
$7$ and $8$. 
The second column displays the number of linear equivalence classes of
subalgebras. The third column has the number of different isomorphism types
of subalgebras. The last column has a field over which all subalgebras
are simultaneously defined.}
\end{center}
\end{table}

Next there is the question of the validation of the results: how can we
be certain that our classifications are correct? Although in this paper
we prove the correctness of the method that we use, there is still ample
possibility to make mistakes while using it.
However we do have some circumstantial evidence for the
correctness of our lists. Firstly, the method does not deal with regular
subalgebras differently than with other subalgebras. But at the end we
find the same regular subalgebras as with Dynkin's algorithm from \cite{dyn}.
Secondly, the $S$-subalgebras that we find in the exceptional types coincide
with the ones found by Dynkin. In the classical types for ranks $\leq 6$ 
we find the same $S$-subalgebras as Lorente and Gruber (\cite{logru})
(except in $D_4$, $D_6$, see below). Thirdly, also the lists of simple 
subalgebras agree with those
found by Dynkin (and in the case of $E_8$ corrected by Minchenko). 

One result of our calculations is that in type $D_{2n}$, for $n=2,3,4$,
there appear 
maximal semisimple subalgebras which are isomorphic, but not linearly
equivalent. In type $D_4$ there are three (linear equivalence classes of) 
maximal subalgebras of types
$A_1B_2$ and $B_3$. In type $D_6$ there are two maximal subalgebras of
each type $A_1C_3$ and $A_5$. And in $D_8$ there are two maximal 
subalgebras of each type $B_2B_2$, $B_4$, $A_1C_4$ and $A_7$. This appears
not to have been known in the literature, for example \cite{mckpat} lists
one algebra for each of the above types. In all cases the algebras are
conjugate under outer automorphisms. It would be interesting to formulate
and prove a general statement about the maximal subalgebras of the Lie algebra
of type $D_{2n}$. However, this would be beyond the scope of this paper.
We intend to come back to it in a subsequent paper. 

This paper is structured as follows. The next section describes a number 
of concepts and results from the literature that we need. This allows us
at the end of the same section
to give a more or less detailed description of the method we use.
The subsequent sections then describe every step in detail.

{\bf Acknowledgement:} It is my pleasure to thank Luigi Scorzato and
Lorenzo Fortunato for having suggested the topic of the paper to me, and
for many fruitful discussions about it.

\section{Preliminaries}

Throughout $\g$ will be a semisimple Lie algebra over $\C$, with a 
fixed Cartan subalgebra $\hh$. The inner automorphism group of $\g$ 
will be denoted $G$. 

\subsection{The action of the Weyl group}\label{p:1}

The Killing form on $\g$ will be denoted $\kappa$; it is defined by
$\kappa(x,y) = \Tr(\ad x\cdot \ad y)$.
The form $\kappa$ is nondegenerate on $\g$ and on $\hh$. Hence we can 
define
a bijection $\hh^*\to \hh$, $\mu\mapsto \hat{\mu}$, by $\kappa(\hat{\mu},h)
=\mu(h)$, for $h\in \hh$. Then $(\mu,\lambda) = \kappa(\hat{\mu},\hat{\lambda})$
defines a non-degenerate bilinear form on $\hh^*$. Also for $\mu\neq 0$ we set
$$\mu^\vee = \frac{2\hat{\mu}}{(\mu,\mu)}.$$
For $\alpha\in \hh^*$ we set $\g_\alpha = \{ x\in \g \mid [h,x]=\alpha(h)x
\text{ for all } h\in\hh\}$.
We let $\Phi$ be the set of all nonzero $\alpha\in \hh^*$ with $\g_\alpha
\neq 0$. Let $\hh_\R^*$ be the real vector space spanned by $\Phi$. Then
$(~,~)$ is an inner product in $\hh_\R^*$, and $\Phi$ is a (reduced) 
root system in $\hh_\R^*$. 

For $\alpha,\beta\in \hh_\R^*$ we set 
$$\langle \alpha,\beta^\vee\rangle = \frac{2(\alpha,\beta)}{(\beta,\beta)}.$$
For $\alpha\in \Phi$ we define the reflection 
$s_\alpha : \hh_\R^*\to \hh_\R^*$ by 
$s_\alpha(\mu) = \mu -\langle \mu,\alpha^\vee\rangle \alpha$. The group
generated by all $s_\alpha$ for $\alpha\in\Phi$ is called the Weyl group,
and denoted $W$. 

Let $\hh_\R$ be the real vector space spanned by 
all $\alpha^\vee$, for $\alpha\in\Phi$.
For $\alpha\in \Phi$ we define the linear map $s_\alpha : \hh_\R
\to \hh_\R$ by $s_\alpha(h) = h-\alpha(h)h_\alpha$. Then
$s_\alpha(\beta^\vee) = \beta^\vee -\langle \alpha,\beta^\vee \rangle 
\alpha^\vee$.
A small calculation shows that the following diagram commutes
\begin{displaymath}
\xymatrix{
\hh_\R^* \ar[r]^{s_\alpha} \ar[d]_{\widehat{\phantom{x}}} & \hh_\R^* 
\ar[d]^{\widehat{\phantom{x}}} \\
\hh_\R \ar[r]^{s_\alpha} & \hh_\R.
}
\end{displaymath}

So, more generally, for $w\in W$ and $\mu \in \hh_\R^*$ we have 
$w\hat \mu = \widehat{w\mu}$.
Also $W$ leaves the Killing form on $\hh_\R$ and on $\hh_\R^*$ invariant.

Let $<$ be an order on $\hh_\R^*$ with 
\begin{itemize}
\item $u < v$ implies $u+w < v+w$ for all $w\in \hh_\R^*$,
\item $u >0$ implies $\lambda u >0$ for all positive $\lambda\in\R$,
and $\lambda u < 0$ for all negative $\lambda\in\R$.
\end{itemize}
We call such a $<$ a {\em root-order}.
A root-order defines a partition $\Phi= \Phi^+\cup \Phi^-$ 
of $\Phi$ into positive and negative roots, and a set $\Delta$ of
simple roots. Conversely, if $\Delta$ is a set of simple roots, then 
we can define an ordering as follows: express $u,v$ as linear combinations
of the elements of $\Delta$, and set $u < v$ if the first nonzero coordinate
of $v-u$ is positive. This order will then yield $\Delta$ as set of
simple roots.

Let $C^*\subset \hh_\R^*$ be set of all $\mu$ with $\langle \mu, \alpha^\vee \rangle
\geq 0$ for all $\alpha >0$.
Then every $W$-orbit in $\hh_\R^*$ has a unique point in $C^*$. It is called
the fundamental Weyl chamber of $\hh_\R^*$.
Also we let $C\subset \hh_\R$ be the set of all $h$ with 
$\alpha(h)\geq 0$ for all $\alpha>0$. Again, every $W$-orbit in $\hh_\R$
has a unique point in $C$. It is called the fundamental 
Weyl chamber in $\hh_\R$.

\subsection{Nilpotent orbits}\label{p:2}
Let $e\in \g$ be nilpotent; then the orbit $G\cdot e$ is called a nilpotent 
orbit. Here we
recall some facts on the classification of nilpotent orbits from
\cite{cart}, \cite{colmcgov}.

Let $e\in \g$ be nilpotent, then by the Jacobson-Morozov lemma there are
$h,f\in \g$ with $[h,e]=2e$, $[h,f]=-2f$, $[e,f]=h$. We say that $(h,e,f)$
is an $\ssl_2$-triple. Note that $G$ acts on $\ssl_2$-triples by
$\sigma\cdot(h,e,f) = (\sigma\cdot h, \sigma\cdot e, \sigma\cdot f)$.
Let $e,e'\in\g$ be nilpotent, lying in $\ssl_2$-triples
$(h,e,f)$ and $(h',e',f')$. Then the following are equivalent:
\begin{itemize}
\item $e,e'$ lie in the same $G$-orbit,
\item $(h,e,f)$ and $(h',e',f')$ lie in the same $G$-orbit,
\item $h,h'$ lie in the same $G$-orbit.
\end{itemize}
Let $(h,e,f)$ be an $\ssl_2$-triple in $\g$. Then $h$ lies in a Cartan
subalgebra of $\g$. As all Cartan subalgebras of $\g$ are $G$-conjugate, after
possibly replacing the triple by a $G$-conjugate, we may assume that $h\in \hh$.
Then $h\in \hh_\R$ (indeed: $\alpha(h)\in \Z$ for all $\alpha\in\Phi$). 
Two elements of $\hh_\R$ are $G$-conjugate
if and only if they are $W$-conjugate (cf. \cite{colmcgov} Theorem 2.2.4). 
Hence, after a further conjugation
we may assume $h\in C$. In fact, this $h$ determines the nilpotent
orbit uniquely; it is called the characteristic of the orbit.

We call an $h\in \hh$ {\em admissible} if it lies in an $\ssl_2$-triple
$(h,e,f)$. Let $e_1,\ldots,e_t$ be representatives of the nilpotent $G$-orbits
in $\g$, lying in $\ssl_2$-triples $(h_i,e_i,f_i)$, with $h_i\in C$.
Then
$$\mathcal{H} = \bigcup_{i=1}^t W\cdot h_i$$
is the set of all admissible elements in $\hh$.

We will often have the need to run through a $W$-orbit $W\cdot h_i$. 
For this Snow (\cite{snow1}, see also \cite{gra6}) has devised an 
efficient algorithm, which makes it possible to run through the orbit and
inspect each element without storing all of the orbit. This feature will
be very important for us.

\subsection{The Dynkin index}\label{p:3}

Assume that $\g$ is simple. It is well-known that 
upto multilication by nonzero scalars, 
there exists a unique nondegenerate symmetric $G$-invariant bilinear form
on $\g$. The Killing form is such a form.

Let $\tilde{\g}\subset \g$ be a simple
subalgebra. Let $\widetilde{G}$ be
the group of inner automorphisms of $\tilde{\g}$. Then $\widetilde{G}
\subset G$. Hence the Killing form $\kappa$ 
of $\g$ induces a $\widetilde{G}$-invariant
bilinear form on $\tilde{\g}$. Let $\tilde{\kappa}$ denote the Killing form
of $\tilde{\g}$. So $\tilde{\kappa}(x,y) =
\eta \kappa(x,y)$ for all $x,y\in \tilde{\g}$, where $\eta$ is a nonzero
scalar.

If we normalise $\kappa$ so that $\kappa(\alpha^\vee,\alpha^\vee)=2$ for
the short roots $\alpha$, and do the same for $\tilde\kappa$, then $\eta$
is called the Dynkin index of $\tilde\g$ in $\g$. It is the same for all
$G$-conjugates of $\tilde\g$. However, it can also happen that nonconjugate
subalgebras have the same Dynkin index.

\begin{lem}\label{lem:1}
Let $\tilde\g\subset \g$ be a semisimple subalgebra that is the direct
sum of simple ideals, $\tilde \g = \tilde \g_1\oplus\cdots\oplus \tilde \g_m$.
Then $\kappa(\tilde \g_i,\tilde \g_j)=0$ for $i\neq j$.
\end{lem}

\begin{pf}
Let $x,y\in \g_i$ and $z\in \g_j$. Then $\kappa([x,y],z) = \kappa(x,[y,z])=0$.
So since $\g_i = [\g_i,\g_i]$ the result follows.
\end{pf}

\subsection{Canonical generators}

Let $\tilde{\g}$ be a semisimple Lie algebra. Then 
$\tilde\g$ has a {\em canonical set of generators}
(\cite{jac}, Chapter IV).
That is a a set of elements $\tilde x_1,\ldots,\tilde x_r$,
$\tilde y_1,\ldots,\tilde y_r$, $\tilde h_1,\ldots,\tilde h_r$ such that
\begin{align}
[\tilde h_i,\tilde h_j] &=0\nonumber\\
[\tilde x_i,\tilde y_j] &= \delta_{ij} \tilde h_i \label{p:4}\\
[\tilde h_j,\tilde x_i] &= \widetilde{C}(i,j)\tilde x_i\nonumber \\
[\tilde h_j,\tilde y_i] &= -\widetilde{C}(i,j) \tilde y_i.\nonumber
\end{align}
Here $\widetilde{C}$ is the Cartan matrix of the root system of $\tilde{\g}$.
We call the sequence $(\tilde h_1,\ldots,\tilde h_r)$ the $h$-part of the
canonical generating set. We note that $\tilde h_i =\alpha_i^\vee$, where
$\{\alpha_1,\ldots,\alpha_r\}$ is a set of simple roots of the root
system of $\tilde\g$. 

Suppose now that $\tilde h_i\in \hh$. Then
the $\tilde h_i\in \hh$ are admissible; hence lie in 
the set $\mathcal{H}$ of Section \ref{p:2}. In the sequel we will say that
the $h$-part of a canonical generating set lies in $\hh$ to mean that all
of its elements do. The next theorem is essentially the same as
\cite{dyn}, Theorem 4.1.

\begin{thm}\label{thm:2}
Let $\widetilde{C}$ be the Cartan matrix of a root system. Let $\tilde x_i$,
$\tilde y_i$, $\tilde h_i$ be elements of a finite dimensional Lie algebra
satisfying the relations (\ref{p:4}). Then the subalgebra generated by
these elements is semisimple, and its root system has Cartan matrix equal
to $\widetilde{C}$.
\end{thm}

\begin{pf}
For $i\neq j$ consider the element 
$$y_{i,j}= (\ad \tilde y_i)^{-\widetilde{C}(j,i)+1}(\tilde y_j).$$
A short calculation (cf. \cite{gra6}, Lemma 7.11.3) shows that 
$[\tilde x_i,y_{i,j}]=0$ and $[\tilde h_i,y_{i,j}] =
(\widetilde{C}(j,i)-2)y_{i,j}$. But $\widetilde{C}(j,i)-2<0$. It follows that
$y_{i,j}$ generates a finite-dimensional irreducible $\ssl_2$-module of
negative highest weight. This is impossible, hence $y_{i,j}=0$. Similarly
we have
$$(\ad \tilde x_i)^{-\widetilde{C}(j,i)+1}(\tilde x_j)=0.$$
Hence the $\tilde x_i$,  $\tilde y_i$, $\tilde h_i$ satisfy the Serre 
relations (see \cite{ser2}, Chapter VI, \S 4). 
This implies that the algebra they generate is a quotient
of the semisimple Lie algebra $\mf{u}$ corresponding to the Cartan matrix 
$\widetilde{C}$ by an ideal. This ideal is the sum of some of the simple
ideals of $\mf{u}$.
But since the $\tilde x_i$,  $\tilde y_i$, $\tilde h_i$ are nonzero, this
ideal has to be zero.
\end{pf}

\subsection{Solving polynomial equations}\label{sec:gb}

In order to construct the subalgebras that we are after, on some occasions
we need to solve polynomial equations in several variables (see Section
\ref{p:8}). For this no general algorithm exists, so we have to do it by
hand. However, a computational tool that makes this a lot easier is provided
by Gr\"obner bases. 

Let $f_1,\ldots,f_s\in k[x_1,\ldots,x_n]$, where $k$ is a field, generate
the ideal $I$. Then solving $f_1=\cdots=f_s=0$ is the same as solving
$g=0$ for all $g\in \mathcal{G}$, where $\mathcal{G}$ is any other generating 
set of $I$.
A Gr\"obner basis is, on many occasions, a particularly convenient
generating set for this purpose. Especially if the Gr\"obner basis 
$\mathcal{G}$ is
computed relative to a lexicographical ordering, then $\mathcal{G}$ 
has a triangular
structure, which often makes solving the polynomial equations easier. Also,
if there are no solutions over the algebraic closure of $k$, then the
reduced Gr\"obner basis is $\{1\}$. So this situation is immediately
detected. Here we do not go into the details, but refer to \cite{clo} for an 
in-depth discussion of Gr\"obner bases and polynomial system solving. 

\subsection{Outline of the method}

Here we summarise the method we use to classify semisimple subalgebras of $\g$.

Let $\widetilde{C}$ be the $r\times r$ Cartan matrix of the root system of a
semisimple Lie algebra. The objective is to classify the semisimple subalgebras
of $\g$ having a root system with Cartan matrix $\widetilde{C}$, up to linear
equivalence. We assume that the semisimple subalgebras of $\g$ of smaller
rank have been classified. We note that the classification for rank $1$
is known from the classification of the nilpotent orbits in $\g$. 

Let $\tilde\g\subset \g$ be a subalgebra with Cartan matrix $\widetilde{C}$,
and canonical set of generators $\tilde h_i$, $\tilde x_i$, $\tilde y_i$,
$1\leq i\leq r$, satisfying (\ref{p:4}). Then the $\tilde h_i$ lie in a Cartan
subalgebra of $\g$. So since all Cartan subalgebras of $\g$ are conjugate
under $G$, we get that $\tilde \g$ is equivalent, and hence linearly
equivalent, to a subalgebra with a canonical generating set with the
$h$-part lying in $\hh$. So we may assume that $\tilde h_i\in \hh$,
and hence $\tilde h_i\in \mathcal{H}$.

In Section \ref{p:11} we describe methods to assemble a set $H$ of 
$r$-tuples $(\tilde h_1,\ldots, \tilde h_r)\in \mathcal{H}^r$ such that
all classes of linearly equivalent subalgebras with Cartan matrix 
$\widetilde{C}$ have a representative that has a canonical generating
set with $h$-part in $H$. Here one of the objectives is to keep this
set ``small''. 

Let $(\tilde h_1,\ldots, \tilde h_r)\in H$.
Section \ref{p:8} contains methods that construct $\tilde x_i$, $\tilde y_i$
in $\g$ satisfying (\ref{p:4}), or decide that no such elements exist.
In the former case we have found a semisimple subalgebra of $\g$ with
Cartan matrix $\widetilde{C}$ by Theorem \ref{thm:2}. In the latter case
the $\tilde h_i$ do not form the $h$-part of a canonical generating set
of a subalgebra with Cartan matrix $\widetilde{C}$.

In Section \ref{p:6} we describe a method for deciding whether two 
semisimple subalgebras are linearly equivalent. So we can get rid of any
linearly equivalent pairs of subalgebras constructed in the previous step.
In fact, linear equivalence depends only on the $h$-parts of the canonical
generating sets; so we can construct the set $H$ so that no linearly
equivalent subalgebras arise. This is important as constructing the
subalgebras is one of the most difficult steps.

We used these methods for classifying the semisimple subalgebras of
the simple Lie algebras of ranks up to $8$. For classifying the 
semisimple subalgebras of the semisimple, but not simple, 
Lie algebras we have a separate
method, described in Section \ref{sec:sss}. Finally the last section 
has the algorithm that we use for deciding inclusion.

\section{Deciding linear equivalence}\label{p:6}

The purpose of this section is to describe an algorithm for deciding whether
two semisimple subalgebras of $\g$ are linearly
equivalent. For this we assume that they are given by canonical sets of
generators, with the $h$-parts lying in $\hh$. First we prove a theorem
that in essence is due to Dynkin (\cite{dyn}, Theorem 1.5). 
Here we show how Dynkin's argument can be adapted to prove the statement
that we need (Corollary \ref{cor:1}).
For this the language of embeddings is more appropriate.

Let $\varphi : \tilde\g\to \g$ be an embedding of the semisimple Lie 
algebra $\tilde\g$ into $\g$. Let $\tilde\hh$ be
a fixed Cartan subalgebra of $\tilde\g$ and assume $\varphi(\tilde\hh) \subset
\hh$. Let $\tilde{\hh}^*_\R$, $\hh_\R^*$ be the $\R$-span of the 
roots of $\tilde\g$ and $\g$ respectively. We define a map 
$\varphi^* : \hh_\R^*\to \tilde\hh_\R^*$ by 
$\varphi^*(\mu)(\tilde h) =
\mu(\varphi(\tilde h))$, where $\tilde h\in \tilde{\hh}^*_\R$. 

Let $\rho : \g \to \gl(U)$ be a representation and
let $\mu$ be a weight of $\rho$, with weight vector $v$. Then for 
$\tilde h\in \tilde{\hh}_\R^*$ we get $\rho(\varphi(\tilde h)) v = \mu
(\varphi (\tilde h) )v = \varphi^*(\mu)(\tilde h)v$. It follows that 
$\varphi^*(\mu)$ is a weight of the representation $\rho\varphi$ of 
$\tilde \g$. In particular, it lies in $\tilde{\hh}_\R^*$. Since the weights
span the spaces $\tilde{\hh}_\R^*$, $\hh_\R^*$, it follows that 
$\varphi^*(\hh_\R^*) = 
\tilde{\hh}_\R^*$.

\begin{lem}\label{p:7}
Let the notation be as above. Fix a root-oder $\prec$ of 
$\tilde{\hh}_\R^*$. Fix also a division $\Phi = \Phi^+\cup \Phi^-$ of
the roots $\Phi$ of $\g$ into positive and negative roots, corresponding to
a root-order of $\hh_\R^*$. Let $\Delta\subset \Phi$ 
denote the corresponding set of simple
roots. Then there exists a root-order $<$ of $\hh_\R^*$, and a $\sigma
\in N_G(\hh)$ with the following properties:
\begin{itemize}
\item The set of positive roots with respect to $<$ is also $\Phi^+$;
\item for $\psi = \sigma\varphi $ we have that $\psi^*(\mu) \prec 
\psi^*(\lambda)$ implies $\mu < \lambda$,
\item if $\psi^*(\mu)\neq \psi^*(\lambda)$ then $\mu <\lambda$ implies
$\psi^*(\mu) \prec \psi^*(\lambda)$. 
\end{itemize}
\end{lem}

\begin{pf}
Let $<'$ be any root-order on $\hh_\R^*$. Define the root-order $<''$ on
$\hh_\R^*$ by 
$\mu <'' \lambda$ if $\varphi^*(\mu) \prec \varphi^*(\lambda)$, or if
those two are equal, $\mu <' \lambda$. Let $\Delta''$ be the corresponding 
set of simple roots. Then there is $w\in W$ with $w\Delta'' = \Delta$.
Let $\sigma \in N_G(\hh)$ be such that the restriction of $\sigma$ to
$\hh$ is $w$. Set $\psi = \sigma\varphi$, and define the root-order $<$ by:
$\mu <\lambda$ if $w^{-1}\mu <'' w^{-1}\lambda$.

Let $\alpha\in \Delta$, and write $\alpha = w\beta$ for some $\beta\in 
\Delta''$. Then $w^{-1}\alpha = \beta >'' 0$. Hence $\alpha >0$; and therefore
the set of positive roots with respect to $<$ is $\Phi^+$.

Next, using $(w^{-1}\mu) (h) = \mu( w h)$ for $h\in\hh$ we get
$\psi^*(\mu) = \varphi^*(w^{-1}\mu)$.
Hence $\psi^*(\mu)\prec \psi^*(\lambda)$ is the same as 
$\varphi^*(w^{-1}\mu) \prec \varphi^*(w^{-1}\lambda)$. This implies that
$w^{-1}\mu <'' w^{-1}\lambda$, and hence $\mu < \lambda$.
The last statement follows directly from the second.
\end{pf}

\begin{thm}\label{thm:3}
Let $\varphi_1,\varphi_1 : \tilde\g\to \g$ be two embeddings of $\tilde\g$
into $\g$. Let $\tilde x_i,\tilde y_i,\tilde h_i$ 
for $1\leq i\leq r$ form a canonical set of generators of $\tilde\g$.
Write $\tilde x_i^1,\tilde y_i^1,\tilde h_i^1$ and 
$\tilde x_i^2,\tilde y_i^2,\tilde h_i^2$ for their images under $\varphi_1$,
$\varphi_2$ respectively. Assume that $\tilde h_i^1,\tilde h_i^2\in \hh$.
Then $\varphi_1$ and $\varphi_2$ are linearly equivalent if and only if
there is a $w\in W$ with $w(\tilde h_i^1) = \tilde h_i^2$ for $1\leq i\leq r$. 
\end{thm}

\begin{pf}
First suppose that $w\in W$ exists.
Let $\rho : \g \to \gl(U)$ be a representation of $\g$.
Let $\mu\in \hh_\R^*$ be a weight of $\rho$, i.e., there are nonzero $u\in U$
with $\rho(h)u = \mu(h) u$ for all $h\in \hh$. Set $\rho_i = \rho\varphi_i$ 
for $i=1,2$. Then $\rho_i$ is a representation of $\tilde\g$.
Observe that $\mu(\tilde h_i^1) = \kappa( \hat\mu, \tilde h_i^1 ) = 
\kappa (w\hat\mu,w(\tilde h_i^1)) = \kappa( \widehat{w\mu}, w(\tilde h_i^1)) 
= (w\mu)(w(\tilde h_i^1))$. But also $w\mu$ is
a weight of $\rho$, with the same multiplicity. Hence it follows that
$\rho_1$ and $\rho_2$ have the same weights with the same multiplicities.
Hence $\varphi_1$, $\varphi_2$ are linearly equivalent.

Now assume that $\varphi_1,\varphi_2$ are linearly equivalent.
By Lemma \ref{p:7} there are $\sigma_1,\sigma_2\in N_G(\hh)$ such that
$\psi_i = \sigma_i\varphi_i$ have the properties stated in Lemma \ref{p:7}
for $\psi$. 

Let $\rho : \g\to \gl(U)$ be an irreducible representation with highest
weight $\lambda$. Set $\rho_i = \rho\circ \psi_i$. Then from Lemma \ref{p:7}
it follows that $\psi_i^*(\lambda)$ is the largest weight of $\rho_i$
in the ordering $\prec$. As the $\psi_i$ are linearly equivalent
it follows that $\psi_1^*(\lambda) = \psi_2^*(\lambda)$. This is the
same as $\lambda( \psi_1(\tilde h) ) = \lambda(\psi_2(\tilde h))$ for
all $\tilde{h} \in \tilde{\hh}$. Now since $\hh_\R^*$ is spanned by dominant
weights, this equality follows for all $\lambda\in \hh_\R^*$. Hence
$\psi_1(\tilde h) = \psi_2(\tilde h)$ for all $\tilde h\in \tilde\hh$.
In particular, this is true for the $\tilde h_i$. So
$\sigma_1 (\varphi_1(\tilde h_i)) = \sigma_2(\varphi_2(\tilde h_i))$.
Now let $w_i\in W$ be such that $\sigma_i |_\hh = w_i$. Then we get
the statement of the theorem with $w = w_2^{-1}w_1$.  
\end{pf}

\begin{cor}\label{cor:1}
Let $\tilde\g_1$, $\tilde\g_2$ be two semisimple subalgebras of $\g$, both
isomorphic to $\tilde\g$. 
Let $\tilde h_1^1,\ldots,\tilde h_r^1$, $\tilde h_1^2,\ldots,\tilde h_r^2$,
be the $h$-parts of canonical sets of generators of $\tilde \g_1$ and
$\tilde \g_2$ respectively. Assume that $\tilde h_i^k\in \hh$ for all $i,k$.
Then $\tilde \g_1$, $\tilde \g_2$ are linearly
equivalent if and only if there is a $w\in W$ with 
$$\{ w(\tilde h_i^1) \mid 1\leq i\leq r\} = 
\{ \tilde h_i^2 \mid 1\leq i\leq r\}.$$ 
\end{cor}

\begin{pf}
The ``if''-part follows immediately from Theorem \ref{thm:3}. For the 
``only if''-part suppose that $\tilde \g_1$, $\tilde \g_2$ are linearly
equivalent. Let $\rho : \g\to \gl(U)$ be a faithful representation.
Then there is $a\in \GL(U)$ with $a \rho(\tilde\g_1)a^{-1} = \rho(\tilde\g_2)$.
Set $\tilde h_i^3 = \rho^{-1}( a\rho(\tilde h_i^1)a^{-1}) \in \tilde \g_2$.
Let $\widetilde{G}_2\subset G$ be the inner automorphism group of $\tilde\g_2$.
Let $\tilde\hh_2^2$, $\tilde\hh_2^3$ denote the subspaces of $\tilde\g_2$
spanned respectively by the $\tilde h_i^2$ and the $\tilde h_i^3$. These
are Cartan subalgebras of $\tilde \g_2$ so there is a 
$\sigma\in \widetilde{G}_2$ with $\sigma(\tilde \hh_2^3) = \tilde \hh_2^2$.
Set $\tilde h_i^4 = \sigma(\tilde h_i^3)$. Then also the $\tilde h_i^4$ form
the $h$-part of a canonical set of generators of $\tilde\g_2$, lying in the
same Cartan subalgebra of $\tilde \g_2$ as the $\tilde h_i^2$. Let $W_2$
denote the Weyl group of $\tilde\g_2$ with respect to $\tilde\hh_2$. Since
different sets of simple roots of $\tilde \g_2$ are conjugate under $W_2$,
there is a $u\in W_2$ such that $\{ u(\tilde h_i^4) \} = \{ \tilde h_i^2\}$.
Let $\tau\in \widetilde{G}_2$ be such that $\tau$ restricted to $\tilde\hh_2$
is $u$. Let $\tilde h_1,\ldots,\tilde h_r$ form the $h$-part of a canonical
generating set of $\tilde \g$.
Let $\varphi_1 : \tilde \g \to \tilde \g_1$ be the isomorphism
sending $\tilde h_i$ to $\tilde h_i^1$. Set $\varphi_2 = \tau\sigma\varphi_1$.
Then $\varphi_2$ is linearly equivalent to $\varphi_1$. Moreover,
$\{\varphi_2(\tilde h_i)\} = \{\tilde h_i^2\}$. Now we get the required
$w\in W$ from Theorem \ref{thm:3}. 
\end{pf}

So we can decide linear equivalence if we can decide whether two sets
of elements, $\{h_1^1,\ldots, h_r^1\}$ and $\{h_1^2,\ldots,h_r^2\}$
of $\hh_\R$ are conjugate under $W$. Since $W$ preserves the Killing form
we assume that the ordering is such that $\kappa(h_i^1,h_j^1) = \kappa
(h_i^2,h_j^2)$ for $1\leq i,j\leq r$. 

We fix a set of positive roots $\Phi^+$ and corresponding set of simple roots
$\{\alpha_1,\ldots,\alpha_l\}$. Then the reflections $s_{\alpha_i}$ generate
$W$.

From Section \ref{p:1} we recall that $C\subset \hh_\R$ is the set of all 
$h$ with $\alpha(h)\geq 0$ for all $\alpha>0$. We note that for a given $h\in
\hh_\R$ it is straightforward to find its unique $W$-conjugate lying in $C$.
Indeed, initially we set $h_0=h$. Let $i\geq 0$ and suppose that $h_i$ is
found. If $h_i\in C$ then we are done. Otherwise there is $\alpha_j$ with
$\alpha_j(h_i)<0$. Then set $h_{i+1} = s_{\alpha_j}(h_i)$. Note that
for $i<j$ we have $h_j-h_i = \sum_{k=1}^l a_k \alpha_k$ with $a_k\in \R$
non-negative, and at least one coefficient $a_k$ is positive. Hence
all $h_i$ are different, and as $W$ is finite the sequence of the $h_i$ 
must land in $C$. From this we also immediately get a $w\in W$ with
$w(h)\in C$.

Next we have a method for deciding whether there is a $w\in W$ with 
$w(h_i^1) = h_i^2$. We first compute $w_1,w_2\in W$ with $w_i(h_1^i)\in C$.
If those are not equal, then the required $w$ does not exist. Otherwise
set $u = w_2^{-1}w_1$; then $u(h_1^1)=h_1^2$. Now the set of all $v\in W$
sending $h_1^1$ to $h_1^2$ is exactly $\Stab_W(h_1^2)u$, where $\Stab_W( h_1^2)$
denotes the stabiliser of $h_1^2$ in $W$.

Set $h= w_2(h_1^2)$. Let $I$ be the set of all $i$ with $\alpha_i(h)=0$.
It is known (cf. \cite{hum3}, Theorem 1.12) that $\Stab_W(h)$ is generated 
by the $s_{\alpha_i}$ with $i\in I$. Now $\Stab_W(h_1^2) = w_2^{-1}\Stab_W(h)
w_2$. This implies that $\Stab_W(h_1^2)$ is generated by the reflections
$s_{w_2^{-1}(\alpha_i)}$, where $i\in I$. The roots $w_2^{-1}(\alpha_i)$ for
$i\in I$ form a simple system of a root subsystem of $\Phi$, of which
$\Stab_W(h_1^2)$ is the Weyl group. 

Now set $h_i^3 = u(h_i^1)$ for $1\leq i\leq r$. We decide if there is
$v\in \Stab_W(h_1^2)$ such that $v(h_i^3) = h_i^2$ for $2\leq i\leq r$.
We can do this as the sequence is shorter. If such a $v$ exists, also
the required $w$ (which is $vu$) exists. In the other case it does not. 

Finally, in order to decide whether there is a $w\in W$ with $\{ w(h_i^1)\}
=\{ h_i^2\}$ we loop over all permutations $\pi$ of $\{1,\ldots,r\}$
with $\kappa( h_{\pi(i)}^1, h_{\pi(j)}^1) = \kappa(h_i^1,h_j^1)$ for
$1\leq i,j\leq r$. For each such $\pi$ we decide whether there is 
a $w\in W$ with $w(h_{\pi(i)}^1) = h_i^2$. Once we find one we stop.

\begin{rem}
This procedure works well in partice if the number of permutations
as above is small. This very often is the case. The main exception
being the case where $(h_1^1,\ldots,h_r^1)$ is the $h$-part of a canonical
generating set of a Lie algebra of type $kA_1$. In situations like that
the algorithm has to work a lot harder, as up to $k!$ permutations have 
to be tried. Fortunately, for the simple Lie algebras
of ranks $\leq 8$ there are not many subalgebras of such a type with large
$k$.
\end{rem}

\begin{rem}
If a class of linearly equivalent subalgebras splits into more than
one class of equivalent subalgebras, then each of the latter classes 
has a representative having a canonical generating set with 
$h$-part that is the same for each of them. Only the other generators
$\tilde x_i$, $\tilde y_i$ differ. 
\end{rem}

\section{Constructing a subalgebra}\label{p:8}

In this section we describe algorithms for constructing a canonical generating
set of a semisimple subalgebra of $\g$, given its Cartan matrix and 
$h$-part.

Let $\tilde{\hh}$ be a subalgebra of $\hh$. For $\mu \in \tilde\hh ^*$ we
set 
$$\g(\mu) = \{ x \in \g \mid [\tilde{h},x] = \mu(\tilde h) x \text{ for all }
\tilde h \in \tilde \hh\}.$$
Then $\g$ is the direct sum of the various $\g(\mu)$.

\begin{lem}\label{lem:orbit}
Let $\mu\in \tilde\hh^*$ be such that $\g(\mu)\neq 0$ and such that there
is a $h\in \tilde\hh$ with $\mu(h)=2$. Set
\begin{align*}
\mathcal{O}_\mu &= \{ u\in \g(\mu) \mid [\g(0),u] = \g(\mu)\},\\
E_\mu  &= \{ e\in \g(\mu) \mid \exists f\in \g(-\mu) \text{ with }
(h,e,f) \text{ is an $\ssl_2$-triple}\}.
\end{align*}
Let $G_0$ be the connected subgroup of $G$ with Lie algebra $\g(0)$.
Then $G_0$ has a dense orbit in $\g(\mu)$, which is equal to $\mathcal{O}_\mu$.
If $E_\mu$ is nonempty then $E_\mu=\mathcal{O}_\mu$.
\end{lem}

\begin{pf}
By standard arguments it is proved that $\kappa$ is non-degenerate
on $\g(-\mu)\oplus \g(\mu)$, and on $\g(0)$. Hence
$$\mf{a} = \bigoplus_{k\in \Z} \g( k\mu )$$
is a reductive $\Z$-graded Lie algebra. In \cite{vinberg} it is shown that
$\g(\mu)$ has a dense $G_0$-orbit. It is clear that a $u\in \g(\mu)$ lies
in this dense orbit if and only if it lies in $\mathcal{O}_\mu$.

Suppose that $E_\mu$ is not empty, and let $e\in E_\mu$. Then from 
$\ssl_2$-representation theory it follows that $\ad e : \g(0)\to \g(\mu)$
is surjective. In other words, $[\g(0),e] = \g(\mu)$. Hence the $G_0$-orbit
of $e$ is dense in $\g(\mu)$. So this last orbit coincides with 
$\mathcal{O}_\mu$. But then also $E=\mathcal{O}_\mu$.
\end{pf}

Let $\widetilde{C}$ be the Cartan matrix of the root system of a semisimple
Lie algebra. Let $\tilde h_1,\ldots,\tilde h_r\in \hh$. We want to find 
$\tilde x_i,\tilde y_i\in \g$ satisfying the relations (\ref{p:4}), or
decide that no such elements exist. We assume that $\tilde h_i\in \mathcal{H}$,
as otherwise the required $x_i,y_i$ certainly do not exist.
The space spanned by $\tilde h_1,\ldots,\tilde h_r$ will be denoted $\tilde\hh$.

First of all, let $\mu_i\in \tilde\hh^*$ be defined by $\mu_i(\tilde h_j) = 
\widetilde{C}(i,j)$. We compute bases of $\g(\mu_i)$ and $\g(-\mu_i)$, and
of $\g(0)$, which is the centralizer of $\tilde\hh$. The $\tilde x_i$,
$\tilde y_i$, if they exist, lie in $\g(\mu_i)$, $\g(-\mu_i)$ respectively.

In the second step we find $\tilde x_1\in\g(\mu_1), \tilde y_1\in\g(-\mu_1)$ 
such that 
$(\tilde h_1,\tilde x_1,\tilde y_1)$ is an $\ssl_2$-triple. For this we use 
Lemma \ref{lem:orbit}. After trying a few random elements
we find an $\tilde x_1\in \g(\mu_1)$ with $[\g(0),\tilde x_1] = \g(\mu_1)$, 
i.e., such that $\tilde x_1$ lies in $\mathcal{O}_{\mu_1}$. 
By solving a set of 
linear equations we either find $\tilde y_1\in \g(-\mu_1)$ such that
$(\tilde h_1,\tilde x_1,\tilde y_1)$ is an $\ssl_2$-triple, or we decide that
that no such $\tilde y_1$ exists. In the latter case there is no 
$\ssl_2$-triple $(\tilde h_1,\tilde x_1,\tilde y_1)$ with $\tilde x_1\in
\g(\mu_1)$, $\tilde y_1\in\g(-\mu_1)$. Indeed, in that case the set
$E_{\mu_1}$ (notation as in Lemma \ref{lem:orbit}) 
is empty. So in the latter case we stop
with the conclusion that the $\tilde x_i,\tilde y_i$ do not exist. 
In the former case we continue.

In this second step we choose a random element $\tilde x_1$. We do stress
that for the existence of the subsequent elements $\tilde x_i,\tilde y_i$, for
$i>1$ it does not matter which $\tilde x_1$ is chosen, as long as 
$[\g(0),\tilde x_1]=\g(\mu_1)$. Indeed: all elements with that property
are conjugate under $G(0)$ as they lie in the same dense orbit.

Now we continue to find the remaining $\tilde x_i$, $\tilde y_i$. 
For this we use two methods, which we call the linear method and the
polynomial method.

For the linear method we suppose that $\tilde x_i,\tilde y_i$, $1\leq i\leq s$,
for a certain $s$ with $1<s <r$, have been found, satisfying (\ref{p:4}). 
We also assume that all different such sets 
are $G$-conjugate. In other words, if $\tilde x_i',\tilde y_i'$ for 
$1\leq i\leq s$ also satisfy (\ref{p:4}), then there exists
$\sigma\in G$ with $\sigma(\tilde x_i') = \tilde x_i$, $\sigma(\tilde y_i') 
= \tilde y_i$,
$\sigma(\tilde h_i) = \tilde h_i$. Note that by the above construction this is 
certainly true for $s=1$.

Set 
$$\g'(\mu_{s+1}) = \{ u\in \g(\mu_{s+1}) \mid [u,\tilde y_i] = 0 \text{ for }
1\leq i\leq s\},$$
$$\g'(-\mu_{s+1}) = \{ u\in \g(-\mu_{s+1}) \mid [u,\tilde x_i] = 0 \text{ for }
1\leq i\leq s\}.$$
Then $\tilde x_{s+1} \in  \g'(\mu_{s+1})$, $\tilde y_{s+1} \in \g'(\-\mu_{s+1})$.
Let also
$\g'(0)$ be the intersection of $\g(0)$ and the centralizer of all 
$\tilde x_i$, $1\leq i\leq s$. By $\ssl_2$-representation theory it follows
that $\g'(0)$ centralises also all $\tilde y_i$, $1\leq i\leq s$. 
Hence $\g'(0)$ acts on $\g'(\mu_{s+1})$. Let 
$G'_0$ be the connected subgroup of $G$ with Lie algebra $\g'(0)$. 
There are now two cases that can occur. 

In the first case, after trying a few random elements, we find a 
$\tilde x_{s+1}\in \g'(\mu_{s+1})$ with $[\g'(0),\tilde x_{s+1}] = 
\g'(\mu_{s+1})$. This means that $G_0'$ has
a dense orbit in $\g'(\mu_{s+1})$. By solving a set of linear equations 
we either find $\tilde y_{s+1}\in \g'(\-\mu_{s+1})$ such that 
$(\tilde h_{s+1},\tilde x_{s+1},\tilde y_{s+1})$ is an $\ssl_2$-triple, or
that no such $\tilde y_{s+1}$ exists. In the former case we say that
the linear method has successfully found $\tilde x_{s+1}$, $\tilde y_{s+1}$.
Note that this also implies that all sets of $\tilde x_i$, $\tilde y_i$
for $1\leq i\leq s+1$ are $G$-conjugate. In the latter case we say that
the linear method has broken down at step $s+1$.

The second case occurs when, after trying a few random elements, we do not
find an $\tilde x_{s+1}$ as above. In this case we also say that the 
linear method has broken down at step $s+1$.

After having found $\tilde x_1$, $\tilde y_1$ we repeat the linear method. 
If it does not break down then in the end we find a complete set of 
$\tilde x_i$, $\tilde y_i$. If it breaks down at step $s+1$, then we 
use the polynomial method.

So for the polynomial method we also assume that $\tilde x_i,\tilde y_i$, for 
$1\leq i\leq s$, have been found, satisfying the relations (\ref{p:4}).
For $s+1\leq k\leq r$ we compute bases of the spaces
$$\g'(\mu_{k}) = \{ u\in \g(\mu_{k}) \mid [u,\tilde y_i] = 0 \text{ for }
1\leq i\leq s\},$$
$$\g'(-\mu_{k}) = \{ u\in \g(\mu_{k}) \mid [u,\tilde x_i] = 0 \text{ for }
1\leq i\leq s\}.$$
We express the $\tilde x_i$, $\tilde y_i$ for $s+1\leq i\leq r$ 
as linear combinations of 
the bases of, respectively, $\g'(\mu_{k})$ and $\g'(-\mu_{k})$, with 
indeterminates as coefficients. Then the $\tilde x_i$, $\tilde y_i$ satisfy
(\ref{p:4}) if and only if certain polynomial equations in the
coefficients are satisfied. We compute the polynomial equations, and 
by Gr\"obner basis techniques (see Section \ref{sec:gb}), 
we either solve them, or decide that no solution exists.

\begin{rem}
Note that the linear method is heuristic in nature. However, it is automatic:
if it succeeds then no further intervention is necessary to construct the 
subalgebra. We note also that there are situations where the linear
method must break down as there are subalgebras that are only defined over
an algebraic extension of $\Q$. In this case using the polynomial method is 
necessary. However, this last method is not entirely automatic (cf. Section
\ref{sec:gb}).
\end{rem}

\begin{rem}
In the next section we give methods to construct a suitable set of candidates
$(\tilde h_1,\ldots,\tilde h_r)$ for the $h$-parts of canonical generating
sets of semisimple subalgebras. This construction is such that 
$(\tilde h_1,\ldots, \tilde h_{r-1})$ will be the $h$-part of a canonical
generating set of a subalgebra of rank $r-1$. However, the $\tilde x_i$,
$\tilde y_i$ for $1\leq i\leq r-1$ do not necessarily lie in the bigger 
subalgebra, as in the two cases (the algebra of rank $r-1$ and of rank $r$)
the spaces $\g(\mu_i)$ are quite different.
\end{rem}

\section{Finding candidates}\label{p:11}

In this section we deal with the problem of finding a suitable set of
candidates for the $h$-parts of canonical generating sets of semisimple
subalgebras of $\g$, with given Cartan matrix $\widetilde{C}$. For this
we first consider a problem involving characters, whose solution will help
us in making the set of candidates smaller. 

Let $\tilde \g$ be a semisimple Lie algebra with canonical generators
$\tilde x_i$, $\tilde y_i$, $\tilde h_i$, $1\leq i\leq r$ 
satisfying (\ref{p:4}). Let $V$ be
a finite-dimensional $\tilde \g$-module. Then $V$ is spanned by common
eigenvectors of the $\tilde h_i$. Moreover, the eigenvalues of the
$\tilde h_i$ are integers. For an $e = (e_1,\ldots,e_r)\in \Z^r$
we set 
$$V_{e} = \{ v\in V \mid \tilde h_i \cdot v = e_i v \text{ for }
1\leq i\leq r\}.$$
Let $x_1,\ldots,x_r$ be indeterminates and write $x^{e} = x_1^{e_1}\cdots
x_r^{e_r}$. Then the polynomial
$$\sum_{e\in \Z^r} (\dim V_e) x^e$$
is called the character of the $\tilde \g$-module $V$. 

Now fix $i$ with $1\leq i\leq r$. For $m\in \Z$ we set $V_m^i = \{v\in V \mid
\tilde h_i  \cdot v = mv\}$, and define the polynomial
\begin{equation}\label{eqn:cp}
f_i(x_i) = \sum_{m\in \Z} (\dim V_m^i) x_i^m.
\end{equation}
We call the polynomial $f_1+\cdots +f_r$ the {\em character-puzzle} of $V$.
It is clear that from the character of $V$ we can compute its character-puzzle.
More generally we say that a polynomial of the form $f_1(x_1)+\cdots
+f_r(x_r)$ is a character-puzzle. It is clear that a character-puzzle does
not necessarily correspond to a character. If it does we say that it is 
{\em solvable}.
Here we consider the following problem: given a character-puzzle
$f=f_1(x_1)+\cdots +f_r(x_r)$ decide whether it is solvable.

For this we proceed as follows. First we note that $V$ is a direct sum
of simple modules, determined by a highest weight, which is an $e=(e_1,\ldots,
e_r)$ with $e_i\geq 0$. From the character-puzzle we retrieve all non-negative
eigenvalues of the $\tilde h_i$. This gives a finite number of possibilities
for the highest weight of a simple constituent of $V$. For each possible
highest weight we compute the character of the corresponding highest weight
module (cf. \cite{gra6}), and from that its character-puzzle $g$. Then we
subtract, $h=f-g$. Then recursively we establish whether $h$ is solvable.

If at least one $h$ that we so obtain is solvable then $f$ itself is solvable.
Otherwise it is not.

Let $\widetilde{C}$ be the Cartan matrix of (the root system of) a semisimple 
Lie algebra $\tilde\g$ of rank $r$. In this section we describe how we find a 
set $H$ of 
$r$-tuples $(\tilde h_1,\ldots,\tilde h_r)\in \mathcal{H}^r$ such that
every semisimple subalgebra of $\g$ isomorphic to $\tilde\g$
is linearly equivalent to a subalgebra with canonical set of generators
$\tilde x_i$, $\tilde y_i$, $\tilde h_i$ with $(\tilde h_1,\ldots,\tilde h_r)\in
H$. We also want the set to be ``small'' (whatever that means). So,
although the set $\mathcal{H}^r$ would be a solution to the problem, it
is far too big. (For example, if $\g$ is of type $E_8$ then it has 
$2611951200^r$ elements.)

A first reduction is given by Corollary \ref{cor:1}: if there are two $r$-tuples
$(\tilde h_1,\ldots,\tilde h_r)$, $(\tilde h_1',\ldots,\tilde h_r')$ such that
there is a $w\in W$ with $w\{\tilde h_i\} = \{\tilde h_i'\}$, then we can
discard one of them. 

Secondly, let $\widetilde{C}_0$ be the $(r-1)\times (r-1)$-matrix in the
top left corner of $\widetilde{C}$. Then we may assume that we know
a set $H_0$ of $(r-1)$-tuples $(\tilde h_1,\ldots ,\tilde h_{r-1})\in
\mathcal{H}^{r-1}$ such that every semisimple subalgebra of $\g$ with Cartan 
matrix $\widetilde{C}_0$ is linearly equivalent to exactly one subalgebra with 
canonical set of generators $\tilde x_i$, $\tilde y_i$, $\tilde h_i$, 
$1\leq i\leq r-1$, with $(\tilde h_1,\ldots,\tilde h_{r-1})\in H_0$.

Therefore we only put $r$-tuples $(\tilde h_1,\ldots,\tilde h_r)$ 
into the set $H$ 
that have $(\tilde h_1,\ldots,\tilde h_{r-1})\in H_0$. Note that for $r=2$ we
know the set $H_0$ from the classification of the nilpotent orbits in $\g$.

So let $(\tilde h_1,\ldots,\tilde h_{r-1})\in H_0$. We want to extend
this $(r-1)$-tuple with an $\tilde h_r$. If we just take any $\tilde h_r\in
\mathcal{H}$, then the set $H$ gets too big. So we perform further reductions.
For this we distinguish two cases. 

In the first case, in the Dynkin diagram of $\widetilde{C}$, the node
labeled $r$ is not isolated. So it is 
connected with $1$, $2$, or $3$ bonds to a simple component $\Gamma_0$ of the 
Dynkin diagram of $\widetilde{C}_0$. Let $\Gamma$ be the simple component
of the Dynkin diagram of $\widetilde{C}$, containing $\Gamma_0$.
Let $i_{1},\ldots,i_s$ be the labels of $\Gamma$, where $i_s=r$. 
Let $\hat\g$ be a simple Lie algebra
with Dynkin diagram $\Gamma$, set of canonical generators $\hat x_i$,
$\hat y_i$, $\hat h_i$, for $1\leq i\leq s$, and Killing form $\hat\kappa$.
As seen in Section 
\ref{p:3} the matrix $(\kappa(\tilde h_{i_k},\tilde h_{i_l}))$ is a
scalar multiple of the matrix $\widehat{B}=
(\hat \kappa(\hat h_i,\hat h_j))$. Furthermore,
we know the scalar factor $\eta$ from comparing 
$\kappa(\tilde h_{i_1},\tilde h_{i_1})$
and $\hat\kappa (\hat h_1,\hat h_1)$. In particular we know what
$\kappa(\tilde h_r,\tilde h_r)$ has to be; denote this value by $\theta$.

Now let $h_1,\ldots,h_t$ be representatives of the $W$-orbits in $\mathcal{H}$
(see Section \ref{p:2}). 
Note that $\kappa(u,u)=\kappa(h_i,h_i)$ for all $u$ in the 
$W$-orbit of $h_i$. So we enumerate the orbits of those $h_i$
such that $\kappa(h_i,h_i)=\theta$. A $\tilde h_r$ in such an orbit is selected
if the matrix $(\kappa(\tilde h_{i_k},\tilde h_{i_l}))$ is equal to $\theta$
times $\widehat{B}$, and $\kappa(\tilde h_r,\tilde h_i)=0$ for $i$ not in
$\{i_1,\ldots,i_s\}$ (cf. Lemma \ref{lem:1}). 

If the number of bonds is $1$ then we can reduce the work further. 
Suppose that the node labeled $i_s=r$ is connected to the node with label
$i_{s-1}$ in $\Gamma$. Let $\beta_1,\ldots,\beta_s$ be the simple roots
of $\hat \g$. Then $\beta_{s-1}$ and $\beta_s$ are conjugate under 
the Weyl group $\widehat W$ of $\hat \g$.  Also, $\beta_i^\vee = \hat h_i$. 
So from what is said in Section \ref{p:1} it follows that  
$\hat h_{s-1}$ and $\hat h_s$ are conjugate under $\widehat W$. Hence
they are conjugate under $\widehat G$, the inner automorphism group of
$\hat \g$. Now an embedding $\hat\g \hookrightarrow \g$ induces an embedding
$\widehat G\hookrightarrow G$. It follows that $\tilde h_{i_{s-1}}$ and 
$\tilde h_r$ must be conjugate under $G$, which implies that they are
conjugate under $W$. The conclusion is that we can limit our search for
suitable elements $\tilde h_r$ to the $W$-orbit of $\tilde h_{i_{s-1}}$. 

In the second case, in the Dynkin diagram of $\widetilde{C}$, the node
labeled $r$ is isolated. In other words, a subalgebra isomorphic to 
$\tilde \g$ is the
direct sum of a subalgebra $\tilde \g_0$, with Cartan matrix $\widetilde{C}_0$,
and a subalgebra isomorphic to $\ssl_2$. Then by Lemma \ref{lem:1}, we can
restrict to adding the $\tilde h_r$ with $\kappa(\tilde h_i,\tilde h_r)=0$
for $1\leq i\leq r-1$. Also in this case we run through $\mathcal{H}$
by enumerating the $W$-orbits of the $h_i$. 

In both cases we can still encounter $W$-orbits that are too large
to enumerate. For example, in order to construct the subalgebras of type
$A_2$, or of type $2A_1$, with the above procedure, one would have to
run through all orbits; when $\g$ is of type $E_8$ this amounts to examining 
$2611951200$ elements. In order to reduce the work needed we use 
character-puzzles.
Let $V$ be the smallest nonzero $\g$-module. For $(\tilde h_1,\ldots,\tilde 
h_r)\in \mathcal{H}^r$ we compute the corresponding character-puzzle, as 
in (\ref{eqn:cp}), where we view $V$ as a $\tilde \g$-module.
We note that all $\tilde h_r$ in the $W$-orbit of $h_i$ lead to the same
character-puzzle. So we decide if the character-puzzle of corresponding
to $(\tilde h_1,\ldots,\tilde h_{r-1},h_i)$ is solvable (i.e., corresponds
to a character of $\tilde\g$). Only if it is, we enumerate the orbit of 
$h_i$.

\begin{rem}
The procedure using character-puzzles eliminates the largest orbits.
For example, for $\g$ of type $E_8$, there are 11 orbits (out of a possible 
69) that need to be enumerated for constructing the subalgebras of type 
$A_2$; they have sizes 240, 2160, 6720, 17280, 30240, 60480, 69120, 181440, 
241920, 483840, 1814400. We also note that $E_8$ is a difficult case in
two respects: it has by far the largest Weyl group, and the largest 
minimal faithful representation of all simple Lie algebras of ranks 
$\leq 8$. The fact that the minimal faithful module has dimension 248 makes
solving the character puzzles rather hard. However, it is still worth
the wile, as the orbits that are excluded this way are so big. From the
sizes of the orbits that still need to be enumerated we also see the
need for an algorithm, as the one of Snow (\cite{snow1}), 
that does so using little memory.
\end{rem}

\section{Subalgebras of semisimple Lie algebras}\label{sec:sss}

Let $\g = \g_1\oplus \g_2$ be the direct sum of two semisimple ideals.
Let $\hh = \hh_1\oplus \hh_2$ be the corresponding decomposition of the
Cartan subalgebra.
Then the Weyl group $W$ of $\g$ is a direct product $W_1\times W_2$,
where $W_1$ (respectively $W_2$) acts trivially in $\hh_2$ (respectively
$\hh_1$). Let $L_i$ be the set of representatives of the linear equivalence
classes of semisimple subalgebras of $\g_i$. We assume that each element
of $L_i$ has a canonical set of generators with $h$-part lying in $\hh$.

Let $\mf{a}\oplus \mf{b}_1$, $\mf{b}_2\oplus \mf{c}$ be elements of 
$L_1$, $L_2$ respectively, where $\mf{b}_1$, $\mf{b}_2$ are isomorphic.
Let $h_i^1,x_i^1,y_i^1$ for $1\leq i\leq s$, $h_i^1,x_i^1,y_i^1$ for 
$s+1\leq i\leq s+r$, $h_i^2,x_i^2,y_i^2$ for $1\leq i\leq r$, 
$h_i^2,x_i^2,y_i^2$ for $r+1\leq i\leq r+m$ be a canonical generating sets 
of respectively
$\mf{a}$, $\mf{b}_1$, $\mf{b}_2$, $\mf{c}$. We assume that the canonical
generators of $\mf{b}_1$, $\mf{b}_2$ are ``in the same order''; that is,
mapping $x_{s+i}^1\to x_i^2$, $y_{s+i}^1\to y_i^2$, $h_{s+i}^1\to h_i^2$,
for $1\leq i\leq r$ defines
an isomorphism $\mf{b}_1\to \mf{b}_2$. Let $\pi$ be a permutation
of $\{1,\ldots,r\}$ preserving the Cartan matrix of $\mf{b}_2$, or,
equivalently, such that $\kappa_2( h_{\pi(i)}^2, h_{\pi(j)}^2) =
\kappa_2(h_i^2,h_j^2)$ for $1\leq i,j\leq r$, where $\kappa_2$ denotes
the Killing form of $\g_2$. Then also mapping $x_{s+i}^1\to x_{\pi(i)}^2$, 
$y_{s+i}^1\to y_{\pi(i)}^2$, $h_{s+i}^1\to h_{\pi(i)}^2$ defines
an isomorphism $\mf{b}_1\to \mf{b}_2$. Let now $\tilde\g$ be the subalgebra
of $\g$ with canonical generating set
\begin{multline*}
 \{ h_i^1,x_i^1,y_i^1 \mid 1\leq i\leq s\} \cup 
\{ h_{s+i}^1+h_{\pi(i)}^2,x_{s+i}^1+x_{\pi(i)}^2,y_{s+i}^1+y_{\pi(i)}^2 
\mid 1\leq i\leq r\} \cup \\\{ h_i^2,x_i^2,y_i^2 \mid r+1\leq i\leq r+m\}.
\end{multline*}
Let $L$ denote the set of subalgebras of $\g$ that can be constructed this 
way. The next proposition is similar to \cite{dyn}, Theorem 15.1.

\begin{prop}\label{prop:sss}
Every semisimple subalgebra of $\g$ is linearly equivalent to an algebra
in $L$. 
\end{prop}

\begin{pf}
Let $\tilde\g$ be a semisimple subalgebra of $\g$. We may assume that
it has a canonical set of generators with $h$-part lying in $\hh$.

Let $p_i : \g\to \g_i$ denote the projection homomorphism. Then 
$\ker p_1 \cap \ker p_2 = 0$. Let $\mf{a}$, $\mf{c}$ be the sum of the ideals 
of $\tilde\g$ that lie respectively in $\ker p_2$ and in $\ker p_1$. 
Let $\mf{b}$ be the sum of the remaining ideals. Then $\tilde\g = 
\mf{a}\oplus \mf{b}\oplus \mf{c}$, with $\mf{a}\subset \g_1$, and
$\mf{c}\subset \g_2$. Let $h_1^1,\ldots,h_s^1$, $h_{r+1}^2,\ldots,h_{r+m}^2$
denote the $h$-parts of canonical generating sets of $\mf{a}$ and $\mf{c}$
respectively. Then $h_i^1\in \hh_1$, $h_i^2\in \hh_2$.

Let $h_i,x_i,y_i$, $1\leq i\leq r$ be a canonical generating set of $\mf{b}$.
Let $H$ denote the set containing the $h_i^1$, $1\leq i\leq s$, 
$h_i^2$, $r+1\leq i\leq r+m$ and $h_i$, $1\leq i\leq r$. We must show
that there is a $w\in W$ such that $w(H)$ is the $h$-part of a subalgebra
in $L$.

Note that $p_1$ and $p_2$ are injective on $\mf{b}$. Write $h_{s+i}^1 = 
p_1(h_i)$, $h_i^2 = p_2(h_i)$, for $1\leq i\leq r$. Those elements form the
$h$-part of a canonical generating set of a semisimple subalgebra
$\mf{b}_1$, respectively $\mf{b}_2$, of $\g_1$ and $\g_2$. Moreover,
the $\mf{b}_i$ are isomorphic to $\mf{b}$.
In particular, $h_1^1,\ldots,h_{s+r}^1$ form the $h$-part of a canonical 
generating set of the subalgebra $\mf{a}\oplus \mf{b}_1$ of $\g_1$. 
Therefore, after possibly reordering the elements of $H$, there is 
a $w_1\in W_1$ such that the $w_1(h_i^1)$ form the $h$-part of a canonical
generating set of an element of $L_1$. Note that this fixes the ordering of
the $h_i\in H$. We can still reorder the $h_i^2\in H$, where 
$r+1\leq i\leq r+m$. So there is a 
$w_2\in W_2$ such that $w_2(h_1^2),\ldots,w_2(h_r^2),w_2(h_{r+1}^2),\ldots,
w_2(h_{r+m}^2)$ form the $h$-part of an element of $L_2$, up to, possibly,
a permutation $\pi$ of the first $r$ elements. This permutation has to 
leave the Cartan matrix of $\mf{b}_2$ invariant.  
\end{pf}

Proposition \ref{prop:sss} gives an immediate procedure for finding
a set $L$ containing representatives of all linear equivalence classes 
of semisimple subalgebras of $\g$. However, it can still happen that 
different members of $L$ are linearly equivalent. For weeding out linear
equivalent pairs we use the algorithm outlined in Section \ref{p:6}.
We also note that $L$ is the disjoint union of two subsets $L'$, $L''$.
Here $L'$ contains the subalgebras that are the direct sum of an algebra
in $L_1$ and an algebra in $L_2$. And $L''$ has the algebras constructed
as above with $\mf{b}_1,\mf{b}_2\neq 0$.
Among the algebras in $L'$ there are
no linear equivalences. Furthermore, an algebra in $L'$ is never linearly
equivalent to an algebra in $L''$.

\section{Deciding inclusion}

For two semisimple subalgebras $\tilde\g_1,\tilde\g_2\subset \g$ we 
write $\tilde\g_1\to \tilde\g_2$ if $\tilde\g_1$ is linearly
equivalent to a subalgebra of $\tilde\g_2$. (Here linear equivalence is 
defined with respecto to $\g$.) Given $\tilde\g_1,\tilde\g_2$,
with canonical generating sets with $h$-parts in $\hh$,
we decide whether $\tilde\g_1\to \tilde\g_2$ in the following way:
First we let $L$ be the set of representatives of the classes of linear
equivalent subalgebras of $\tilde\g_2$. We get this from the classification
of those subalgebras of $\tilde\g_2$. All are given by canonical generating
sets having $h$-parts in $\hh$. Then we decide whether $\tilde \g_1$ is
linearly equivalent to an element of $L$, using the algorithm from Section
\ref{p:6}.

Now let $\tilde\g_1,\ldots,\tilde\g_s$ be a chain of subalgebras. 
This means that $\tilde \g_i\to \tilde \g_{i+1}$ for 
$1\leq i<s$. Then we can compute a {\em realization} of the chain; that is 
if necessary we replace the $\tilde\g_i$ by linear conjugates such that 
$\tilde \g_i\subset \tilde \g_{i+1}$. For this we start ``at the top'', and 
suppose that $\tilde \g_i \to \cdots \to \tilde \g_s$ has been realised.
We compute the subalgebras of 
$\tilde\g_{i}$ isomorphic to $\tilde\g_{i-1}$, up to linear equivalence.
We find a subalgebra $\mf{s}$ that is linearly equivalent to $\tilde \g_{i-1}$,
and replace $\tilde \g_{i-1}$ by $\mf{s}$.

\def\cprime{$'$} \def\cprime{$'$} \def\Dbar{\leavevmode\lower.6ex\hbox to
  0pt{\hskip-.23ex \accent"16\hss}D} \def\cprime{$'$} \def\cprime{$'$}
  \def\cprime{$'$}

\end{document}